\documentclass [12pt] {article}
\usepackage{amsmath}
\usepackage{amsopn}
\usepackage{amssymb}
\usepackage{latexsym}
\usepackage{amsfonts}
\usepackage{amsthm}
\newtheorem{thm}{Theorem}[section]

\newtheorem{pro}[thm]{Proposition}
\newtheorem{lem}[thm]{Lemma}
\newtheorem{cla}[thm]{Claim}

\newtheorem{cor}[thm]{Corollary}

\theoremstyle{definition}
\newtheorem{obs}[thm]{Observation}

\newtheorem{exa}[thm]{Example}
\newtheorem{defn}[thm]{Definition}

\newtheorem{conj}[thm]{Conjecture}

\newcommand{\een}{\end{enumerate}}
\newcommand{\blem}{\begin{lem}}
\newcommand{\elem}{\end{lem}}
\newcommand{\bcl}{\begin{cla}}
\newcommand{\ecl}{\end{cla}}
\newcommand{\ethm}{\end{thm}}
\newcommand{\bpr}{\begin{pro}}
\newcommand{\epr}{\end{pro}}
\newcommand{\bco}{\begin{cor}}
\newcommand{\eco}{\end{cor}}
\newcommand{\bcon}{\begin{conj}}
\newcommand{\econ}{\end{conj}}
\newcommand{\bde}{\begin{defn}}
\newcommand{\ede}{\end{defn}}
\newcommand{\bex}{\begin{exa}}
\newcommand{\eexa}{\end{exa}}
\newcommand{\bobs}{\begin{obs}}
\newcommand{\eobs}{\end{obs}}
\newcommand{\bexe}{\begin{exe}}
\newcommand{\eexe}{\end{exe}}

\def\bn{\bf n}
\date{}

\title{On the Statistical Independence of Shift-Register Pseudorandom Multisequence over Part of the Period}

\author {Mordechay B. Levin and Irina L. Volinsky}
\begin{document}
\maketitle
\baselineskip=18pt
\hsize=15.0 true cm
\hoffset=-0.8 true cm 

\begin{abstract}
  In this paper we construct a pseudorandom multisequence  $(x_{n_1,...,n_r})$
    based on $k$th-order  linear recurrences modulo $p$,
  such that the discrepancy of the $s$-dimensional  multisequence $(x_{n_1+i_1,...,n_r+i_r})_{1 \leq i_j \leq s_j,  1 \leq j \leq r
  }$  $1 \leq n_j \leq N_j,  1 \leq j \leq r$ is equal to $O((N_1 \cdots N_r)^{-1/2} \ln^{s+3r}(N_1 \cdots N_r))$, where $s=s_1 \cdots s_r$, for all $N_1,...,N_r$ with $1 < N_1 \cdots N_r \leq p^k$.
\end{abstract}

\section{Introduction}
\setcounter{equation}{0}

 Equidistribution and statistical
independence properties of uniform
pseudorandom numbers can be analyzed based
on the discrepancy of certain point sets in
$[0,1)^s$:

Let $\mathbf{x}_n=(x_{n,1}, \ldots, x_{n,s}), \;\;\; n =0, \ldots, N-1$,  be a sequence of points in an $s$-dimensional unit cube $[0, 1)^s$; $v = [0, \gamma_1)\times \ldots \times [0, \gamma_s)$ a box in $[0, 1)^s$.
The quantity:
\[
D((\mathbf{x}_n)_{n=0}^{N-1}) = \sup_{0 < \gamma_1, \ldots, \gamma_s \leq 1} \left| \#\left\{n \in [0, N-1]\;\Big|\; \mathbf{x}_n \in v\right\}/N - \gamma_1 \ldots \gamma_s\right|
\] 
is called the {\it discrepancy} of $(\mathbf{x}_n)_{n=0}^{N-1}$.

 Let us consider pseudorandom numbers (abreviated PRN) obtained by means of the {\it shift-register} method:
 
 Let $p$ be a prime, let $k \geq 2$ be an integer, and generate a $k$th-order linear recurring sequence $y_0,y_1,...\; \in F_p$ by
\begin{equation} \label{1}
 y_{n+k}  \equiv a_{k-1} y_{n+k-1} +...+ a_0y_n \; \mod p,   \qquad n=0,1,...,  
\end{equation}
where $y_0,...,y_{k-1}$ are initial values that are not all zero. The integer coefficients $a_0,...,a_{k-1}$ in (\ref{1}) are 
chosen in such a way that, if they are viewed as elements of the finite field $F_p$, then the  characteristic polynomial
$$
    f(x)=x^k-a_{k-1}x^{k-1}- ... -a_0 \in F_p[x]
$$ 
of the recursion (\ref{1}) is a primitive polynomial over $F_p$. Note that the characteristic polynomial $f$ has a root $\beta$
 in the extension field $F_q$ of $F_p$, where $q=p^k$. 
Let $F_q^* = F_q \setminus \{0\}$ be the multiplicative group of nonzero elements of $F_q$ and let
 $\hat{F}_q =\left \{ \beta \in F^*_q \;\;\;| \;\;\;\beta -\rm{is}\;\; \rm{a}\;\; \rm{primitive}\;\; \rm{root} \right\}$.
We see that $\#F_q^*=q-1$, and
$\#F_q = \varphi(q-1)$, where $\varphi$ is the Euler's function. 
Let $Tr$ denote the trace function from $F_q$ to $F_p$. 
 It is known (see, e.g., [Ni, p. 212]) that there exists an $\alpha \in F_q^{*}$
 such that
\begin{equation} \label{700}
     y_n =Tr(\alpha \beta^n) \qquad {\rm for} \qquad n=0,1,...
\end{equation}

 In the {\it digital multistep method}, the sequence  $y_0,y_1,...$ is transformed into a sequence $x_0,x_1,...$
 of uniform PRN in the following way
\begin{equation}
x_n = \sum_{j=1}^{k}   y_{kn+j-1}/p^{j}, \qquad {\rm for} \qquad  n=0,1, \ldots     \label{lab2}
\end{equation} 

 In a series of papers, Niederreiter (see the review in [Ni]) 
 proved that there exists a  characteristic polynomial $f$ such that 
$$  D ((x_n,...,x_{n+s-1})_{n= 0}^{N-1})  = O(\frac{\sqrt{ \tau}(\ln \tau)^{s+1}}{N}), \quad {\rm for} \quad N=1,...,\tau,  $$
where $\tau$ is the period of the sequence of pseudorandom numbers.\\

This estimate is interesting for $ N \gg \sqrt{ \tau}(\log \tau)^{s+1}$. 
 In [Le1],[Le2], Levin descibed a class of uniform PRN sequences $(z_n)_{n \geq 0}$, having  
 a nontrivial discrepancy estimates also for a small part of the period:
\begin{equation}
 D ((z_n,...,z_{n+s-1})_{n= 0}^{N-1})  = O(N^{-1/2}  \ln^{s+3} N ), \quad {\rm for} \quad N=1,2,...     \label{2}
\end{equation} 

 Our goal is to obtain a nontrivial discrepancy estimate similar to (\ref{2}) for a small part of the period for sequences, some subsequences 
 and multisequences of PRN
   based on $k$th-order  linear recurrences modulo $p$.
The method of the proof is based on Korobov's approach [Ko2] (see also [Le1], [Le2]).
 Similar results can be obtained for the pseudorandom sequences described in [Le2] and [Le3].\\

In this paper we will prove the following theorems: \\

{\bf Theorem 1.} {\it Let  $\epsilon \in (0,1)$, and let $\mathbf{x}(n) = (x_n, \ldots, x_{n+s-1})$ (see (\ref{700}), (\ref{lab2})).
 Then there exist more than $(1-\epsilon) q \varphi(q-1)$ pairs 
 $(\alpha, \beta) \in F_q \times  \hat{F}_q$ such that for any $N \in [1,q]$ and any $b_1,b_2 \in [1,q]$,   the bound
\begin{equation}  \nonumber
  D((\mathbf{x}(b_1n +b_2))_{n=0}^{N-1}) \leq \epsilon^{-1} c_1 N^{-1/2} \ln^{s+2.5}  6N  \ln^{2.5} 3 \ln 6N,
\end{equation}
holds, where the constant $c_1$ depends only on $s,b_1$ and $b_2$.}
\\
\\
\\Let $r \geq 2$, $\boldsymbol{\alpha} = (\alpha_1, \ldots, \alpha_r) \in F_q^r$, and $\boldsymbol{\beta} = (\beta_1, \ldots, \beta_r) \in \hat{F}_q^r$, $\mathbf{n} = (n_1, \ldots, n_r)$ and let
\[
y(\mathbf{n}, j) \equiv  \sum_{i=1}^{r} Tr \left(\alpha_i \beta_i^{k (n_i + s \tilde{n}_i) + j-1} \right)  \mod p, \quad 
          y(\mathbf{n}, j)\in [0,p-1],
\]
where,
$$
    \tilde{n}_i \equiv \sum_{w \in [1, r] ,w \neq i} n_w s^{d_{w,i}}  \mod s^{r-1} ,\;\; \tilde{n}_w \in [0,s^{r-1} ),
$$
with
\begin{equation} \label{800}
          d_{w,i} =w-1 \;\; {\rm for} \;\; w<i,  \quad {\rm and } \quad d_{w,i} =w-2 \;\; {\rm for} \;\; w>i.
\end{equation}
Consider the following multisequence
\begin{equation}
x_\mathbf{n} = \sum_{j=1}^{k} \frac{y(\mathbf{n}, j) }{ p^j},\;\; \rm{and} \;\; \mathbf{x}_{\mathbf{n}} = (x_{\mathbf{n}+\mathbf{i}})_{
 \substack{0 \leq i_w < s_w \\ 1 \leq w \leq r}},   \; {\rm with} \; s_i \in [1,s], 1 \leq w \leq r,  \label{w1}
\end{equation}
where $\mathbf{i} = (i_1, \ldots, i_r)$.\\

{\bf Theorem 2.} {\it
Let $\epsilon \in (0,1)$. Then there exist more than $(1-\epsilon)q^r (\varphi(q-1))^r$ values $(\boldsymbol{\alpha}, \boldsymbol{\beta}) \in F_q^r \times \hat{F}_q^r$ such that, for any $N_i \in [1, q],\; 1 \leq i \leq r$, $N_1 \ldots N_r \leq q$, and any 
  $s_i \in [1, s],\; 1 \leq i \leq r$, $s_0=s_1 \ldots s_r \leq s$, the bound
$$
D ((\mathbf{x}_{\mathbf{n}})_{ \substack{ 0 \leq n_w < N_w  \\ 1 \leq w \leq r}}) 
\leq  \epsilon^{-1}c_2 \left(N_1 \ldots N_r \right)^{-1/2} 
 \ln^{s_0+2.5r} (2^{r+1} N_1 \ldots N_r)     \ln^{2.5r} 3\ln (6N_1 \cdots N_r)
$$
holds, where the constant $c_2$ depends only on $s_1,...,s_r$.}
 
\section{Auxiliary results}
For the integer $b \geq 2$, let denote $C(b) = (-b/2, b/2]\bigcap Z$.
Let $C_s(b)$  be the inner product of $s$ copies of $C(b)$.
Consider point sets for which all coordinates of all points have a finite digit expansion in a fixed base $b \geq 2$. Let
\begin{equation}
\mathbf{w}_n = (w_n^{(1)}, \ldots, w_n^{(s)}) \in [0, 1)^s, \;\;\;\;\; n =0,1,\ldots, N-1,  \label{a1}
\end{equation}
where, for an integer $m \geq 1$, we have
\[
w_n^{(i)} = \sum_{j=1}^{m} w_{nj}^{(i)} b^{-j}, \;\;\;\;\; 0 \leq n \leq N-1,\;\; 1 \leq i \leq s,  
\]
with $w_{nj}^{(i)} \in Z_b$ for $0 \leq n \leq N-1,\;\; 1 \leq i \leq s,\;\; 1 \leq j \leq m$.
For $(h_1, \ldots, h_m) \in C_m(b)$, define
\[
d(h_1, \ldots, h_m) = \begin{cases} {\rm largest}\;\;\; d\;\;\;{\rm with}\;\;\;h_d \neq 0,\;\;\;\;\;{\rm if}\;\;\;(h_1, \ldots, h_m) \neq \mathbf{0},\\
    0, \;\;\;\;\;\;\;\;\;\;{\rm if}\;\;\;(h_1, \ldots, h_m) = \mathbf{0}.
    \end{cases}
\]
For $b = 2$, put $Q_b(h_1, \ldots, h_m) = 2^{-d(h_1, \ldots, h_m)}$.
For $b > 2$, put
\[
Q_b(h_1, \ldots, h_m) = \begin{cases} b^{-d}(\rm{csc} \frac{\pi}{b}|h_d| + \sigma(d,m)),\;\;\;\;\;{\rm if} \;\;\;(h_1, \ldots, h_m) \neq \mathbf{0},\\
    1, \;\;\;\;\;\;\;\;\;\;{\rm if}\;\;\;(h_1, \ldots, h_m) = \mathbf{0},
    \end{cases}
\]
where $d = d(h_1, \ldots, h_m)$ and where $\sigma(d,m) = 1$ for $d < m$ and $\sigma(m,m) = 0$.
\\Let $C(b)^{s \times m}$ be the set of all $s \times m$ matrices with entries in $C(b)$, $C^*(b)^{s \times m} = C(b)^{s \times m} \setminus \{ 0\}$. For $H = (h_{ij}) \in C(b)^{s \times m}$, we define
\[
W_b(H) = \prod_{i=1}^{s} Q_b(h_{i1}, \ldots, h_{im}).
\]
\\\\{\bf Theorem A.} [Ni, Theorem 3.12, p.37, Lemma 4.32, p.68] {\it
If $P$ is the point set (\ref{a1}), and $m \geq [\log_b N]$ then
\[
D (P) \leq \frac{s b}{N} + \sum_{H \in C^*(b)^{s \times m}} W_b(H) \left| \frac{1}{N} \sum_{n=0}^{N-1} e \left(\frac{1}{b} \sum_{i=1}^{s} \sum_{j=1}^{m} h_{ij} w_{nj}^{(i)}\right)\right|,
\]
where } $e(x) = e^{2 \pi \sqrt{-1} x}$.
\\\\{\bf Lemma 1.} [Ni, Lemma 5, p.18]  {\it
Let $s \geq 1$ and $m \geq 1$ be integers. Then} 
\begin{equation}
\sum_{H \in C^*(b)^{s \times m}} W_b(H) < \left(\frac{2}{\pi} m \ln b + \frac{7}{5} m - \frac{m-1}{b} \right)^s.  \label{aa1}
\end{equation}
\\\\{\bf Lemma 2.} (see e.g., [KoSh, p.9, p.13, ref. 3.3]) {\it
Let $\beta \in F_q$,
\[
\delta(\beta) = \begin{cases} 1,\;\;\;\;\;if \;\;\;\beta =0,\\
    0, \;\;\;\;\;\;\;\;\;\;\rm{otherwise}.
    \end{cases}
\]
Then}
\[
\delta(\beta) = \frac{1}{q} \sum_{\alpha \in F_q} e \left(\frac{Tr(\alpha \beta)}{p} \right).
\]
\\\\For proof of the following well known lemma see, e.g., [Ko, p.13], or [LeVo, Lemma 7, p.156].
\\{\bf Lemma 3.}  {\it
Let $N \in [0, T-1],\;\;\; T \in [1, q]$ and $x_n$ be a real ($0 \leq n \leq T-1$). Then
\[
\left| \sum_{n=0}^{N-1} e(x_n) \right|  \leq  \sum_{m= -T/2}^{T/2} \frac{1}{\bar{m}} \left|\sum_{n=1}^T e\left(x_n + \frac{n m}{T}\right)  \right|,
\]
where } $\bar{m} = \max (1, |m|)$.
It is easy to see  that
\begin{equation}
\sum_{m = -T / 2}^{T / 2} \frac{1}{\bar{m}} \leq 3 + 2 \ln T.    \label{aa2}
\end{equation}
\\\\{\bf Lemma 4.} {\it
Let $N_i \in [0, T_i - 1],\;\;\; 1 \leq i \leq r$, where $T_i \in [1, q], \;\;\; 1 \leq i \leq r$, and $\mathbf{\bar{n}} = (n_1, \ldots, n_r)$. Then}
\[
\left| \sum_{n_1=0}^{N_1-1} \ldots \sum_{n_r=0}^{N_r-1} e(x_{\bn}) \right|  \leq  \sum_{m_1 = -T_1/2}^{T_1/2} \ldots \sum_{m_r = -T_r/2}^{T_r/2} \frac{1}{\bar{m}_1 \ldots \bar{m}_r} 
\]
\[
\times \left|\sum_{n_1=0}^{T_1-1} \ldots \sum_{n_r=0}^{T_r-1} e\left(x_{\bn} + \frac{n_1 m_1}{T_1} + \ldots + \frac{n_r m_r}{T_r}\right) \right|.
\]
The proof of the Lemma 4 is the same as the proof of the Lemma 3.
\\\\{\bf Lemma 5.}  {\it
Let $r \geq 2$, $s_1, \ldots, s_r \in [1,s]$ be integers, 
\begin{equation}
\widetilde{n_l + i_l} \equiv \sum_{\substack{ w \in [1, r]  \\w \neq l }} (n_w + i_w)  s^{d_{w,l}}  \mod s^{r-1} ,\;\;\;\;\;\widetilde{n_l + i_l} \in [0, s^{r-1}).
     \label{cd}
\end{equation}
Then for } $l \in [1, r]$,
\[
  \# \left\{k (s \widetilde{n_l + i_l} + i_l) + j \Big| 0 \leq j < k, \; 0 \leq i_{\nu} < s_{\nu},\; \nu = 1, \ldots, r \right\} = ks_1 \dots s_r.
\]
{\bf Proof.}
By (\ref{cd}) it is enough to prove that there are no two vectors $(i_1, \ldots, i_r, j) \neq (i_1^{'}, \ldots, i_r^{'}, j^{'})$ with
\begin{equation}
k (s\widetilde{n_l + i_l} + i_l) + j = k (s \widetilde{n_l + i_l^{'}} + i_l^{'}) + j^{'}. \label{bb1}
\end{equation}
We see $j \equiv j^{'} \mod k $. Hence $j = j^{'}$. By (\ref{bb1}), we get:
\begin{equation}
i_l \equiv i_l^{'} \mod s.   \label{bb2}
\end{equation}
Therefore $i_l = i_l^{'}$. From (\ref{bb1}) we have that $\widetilde{n_l + i_l} = \widetilde{n_l + i_l^{'}}$. Hence
\[
\sum_{\substack{ w \in [1, r]  \\w \neq l }} (n_w + i_w) s^{d_{w,l}}
\equiv \sum_{\substack{ w \in [1, r]  \\w \neq l }} (n_w + i_w^{'})  s^{d_{w,l}} \mod s^{r-1}.
\]
Thus
\begin{equation}
\sum_{\substack{ w \in [1, r]  \\w \neq l }} (i_w - i_w^{'})  s^{d_{w,l}} \mod s^{r-1}.    \label{bb3}
\end{equation}
For $l \neq 1$, we obtain
\begin{equation}
i_1 - i_1^{'} \equiv 0 \mod s.   \label{bb4}
\end{equation}
Note that for $l=1$, (\ref{bb4}) follows from (\ref{bb2}). Hence $i_1 = i_1^{'}$. Suppose that $i_j = i_j^{'}$ for $j=1, \ldots, \nu-1$, 
$\nu \neq l$, $\nu \geq 2$. By (\ref{800}) we get $d_{w,l} -d_{\nu,l} \geq 1$ for $w>\nu$.
 We deduce from (\ref{bb3}) that
\[
i_{\nu} - i_{\nu}^{'} + \sum_{\substack{ w \in (\nu, r]  \\ w \neq l }}  (i_w - i_w^{'} )
  s^{d_{w,l} -d_{\nu,l}}  \equiv 0 \mod s^{r-1 -d_{\nu,l}},
\]
\[
i_{\nu} - i_{\nu}^{'} \equiv 0 \mod s_{\nu}, \;\;\;\;\; {\rm and}\;\;\;\;\; i_{\nu} = i_{\nu}^{'}. 
\]
By induction, Lemma 5 is proved. \qed
\\\\{\bf Lemma 6.}
Let $q = p^k > 3000$, $T \in [1, 4q]$. Then
\begin{equation}
\frac{kT}{\varphi(q-1)} \leq 40 \ln 3T  \ln (3\ln 3T).   \label{i2}
\end{equation} 
{\bf Proof.} First we will prove that
\begin{equation}
\frac{T}{\varphi(q-1)} \leq 10 \ln 3\ln3 T.   \label{i1}
\end{equation} 
By [Sa, p.15, ref. 3a; p.9, ref. 2]
\[
\frac{n}{\varphi(n)} < e^{0.58} \ln \ln n + \frac{2.6}{\ln \ln n},\;\;\;\;\; \rm{and}\;\;\;\;\; \varphi(n) \geq n^{2/3} \;\;\;\;\; n \geq 30.
\]
Hence
\begin{equation}
\frac{q}{\varphi(q-1)} \leq \frac{q-1}{\varphi(q-1)} + \frac{1}{(q-1)^{2/3}} < 1.5 + 2 \ln \ln (q-1). \label{i4}
\end{equation}
If $T < \varphi(q-1)$, then (\ref{i1}) is true. Let $T \geq \varphi(q-1)$, then 
\[
\ln T \geq \frac{2}{3} \ln (q-1),\;\;\;\;\; \rm{and}\;\;\;\;\; \ln \ln T \geq \ln \ln (q-1) - 0.5.
\]
By (\ref{i4}) we have
$$
    \frac{T}{\varphi(q-1)} \leq   4(  1.5 + 2 \ln \ln (q-1)  )  \leq 10 + 8\ln \ln 3T.
$$
The inequality (\ref{i1}) is proved. If $T \leq \varphi(q-1) / k$, then (\ref{i2}) is true. 
Let $T \geq \varphi(q-1)/ k  \geq  \varphi(q-1) / \log_2 q$, then
\[
\ln T \geq \ln q - \ln (1.5 + 2 \ln \ln (q-1)) - \ln (\log_2 q).
\]
Hence
\[
 4\ln T \geq \theta(q), \;\rm{where}\; \theta(x) = 4 \ln x - 4 \ln (1.5 + 2 \ln \ln (x-1)) - 4 \ln (\log_2 x). 
\]
It is easy to verify that $\left(\theta(x) - \log_2 (x)\right)^{'} < 0$ for $x > 3000$, and that $\theta(3000) > 15> \log_2 3000$.
Thus
\[
 4 \ln T > \log_2 q.
\]
Applying (\ref{i1}) we get (\ref{i2}).  Lemma 6 is proved. \qed
\section{Proof of Theorem 1}
By (\ref{lab2}) and Theorem A, with $m = [\log_p N]$ we get:
\begin{equation}
N \cdot D ((\mathbf{x}(b_1n+b_2))_{n=0}^{N-1}) \leq s p 
+ \sum_{H \in C^*(b)^{s \times m}} W_p(H) \Big| S(H) \Big|,   \label{lab3}
\end{equation}
where
\[
S(H) = \sum_{n=0}^{N-1} e\left(\frac{1}{p} \sum_{i=0}^{s-1} \sum_{j=1}^{m} h_{ij} Tr\left(\alpha \beta^{k((b_1n+i)+b_2) +j-1}\right)\right). 
\]
Using Lemma 3, we have
\[
\Big| S(H) \Big| \leq \sum_{m_1 = -T/2}^{T/2} \frac{1}{\bar{m}_1} \Big|\dot{S}(H, m_1, \alpha, \beta) \Big|,
\]
where
\[
\dot{S}(H, m_1, \alpha, \beta) = \sum_{n=0}^{T-1}  e\left(\frac{1}{p} Tr\left(\alpha \beta^{k b_1n} \gamma(\beta,H)\right) + \frac{m_1 n}{T}\right)
,\]
with
\begin{equation}
\gamma(\beta, H) = \sum_{i=0}^{s-1} \sum_{j=1}^{m}  h_{ij} \beta^{kb_1 i + kb_2+ j -1},\;\;\;\;\;\rm{and} \;\;\;\;\;T \in [N, q].   \label{lab4}
\end{equation}
Taking $m = [\log_p T]$, we obtain from (\ref{lab3}) and (\ref{lab4}), that 
\begin{equation}
N \cdot D ((\mathbf{x}(b_1n+b_2))_{n=0}^{N-1}) \leq s p + T \tilde{D}_T(\alpha, \beta), \label{aa}
\end{equation}
where,
\begin{equation}
T \tilde{D}_T(\alpha, \beta) = \sum_{H \in C^*(b)^{s \times m}} W_p(H) \sum_{m_1 = - T/2}^{T/2} \frac{1}{\bar{m}_1}  \left| \dot{S}(H,m,\alpha, \beta) \right|.   \label{5a}
\end{equation}
Let,
\begin{equation}
\chi = \frac{1}{q \cdot \varphi(q-1)} \sum_{\beta \in \hat{F}_q} \sum_{\alpha \in F_q}  \left| \dot{S}(H,m,\alpha, \beta) \right|.   \label{5c} 
\end{equation}
\\Using the Cauchy - Shwartz inequality, we get:
\[
\chi^2 \leq \frac{1}{q \cdot \varphi(q-1)} \sum_{\beta \in \hat{F}_q} \sum_{\alpha \in F_q} \left|\sum_{n = 0}^{T-1} e\left(\frac{1}{p}  \cdot Tr\left(\alpha \beta^{k b_1 n}\gamma(\beta,H)\right) + \frac{m_1 n}{T}\right) \right|^2 
\]
\[
= \frac{1}{q  \varphi(q-1)} \sum_{\beta \in \hat{F}_q} \sum_{\alpha \in F_q}  \sum_{n_1, n_2 = 0}^{T-1} e\left(\frac{Tr\left(\alpha \beta^{kb_1n_1} \gamma(\beta,H) \right) - Tr\left(\alpha \beta^{kb_1 n_2} \gamma(\beta,H) \right)}{p}\right)
\]
\[
\times e\left(\frac{m_1 (n_1 - n_2)}{T}\right).
\]
Applying Lemma 2, we obtain:
\[
\chi^2 = \frac{1}{\varphi(q-1)} \sum_{\beta \in \hat{F}_q}\sum_{n_1, n_2 = 0}^{T-1} \delta \Big(\gamma(\beta,H) \left(\beta^{k b_1n_1} - 
  \beta^{k b_1 n_2}\right) \Big)
e \left(\frac{m_1 (n_1 - n_2)}{T} \right)
\]
\[
\leq \frac{1}{\varphi(q-1)} \sum_{\beta \in \hat{F}_q}\sum_{n_1, n_2 = 0}^{T-1} \delta \Big(\gamma(\beta,H) \left(\beta^{kb_1 n_1} - \beta^{k b_1 n_2}\right) \Big).
\]
Bearing in mind that $\delta(\gamma_1 \gamma_2) \leq \delta(\gamma_1) + \delta(\gamma_2)$, we have:
\begin{equation}
\chi^2 \leq \chi_1 + \chi_2,    \label{aa16}
\end{equation}
where
\[
\chi_1 = \frac{1}{\varphi(q-1)} \sum_{\beta \in \hat{F}_q}\sum_{n_1, n_2 = 0}^{T-1} \delta \Big(\beta^{k b_1n_1} - \beta^{k b_1n_2} \Big),\;\;\;\;\; \rm{and} \;\;\;\;\;
\chi_2 = \frac{T^2 \xi}{\varphi(q-1)},
\]
with
\begin{equation}
\xi = \sum_{\beta \in \hat{F}_q} \delta \Big(\gamma(\beta,H) \Big).   \label{aa7}
\end{equation}
Consider $\chi_1$. Let $\beta^{k b_1n_1} - \beta^{k b_1n_2}  =0$. Then $\beta^{kb_1(n_1 - n_2)} = 1$. Taking into account that $\beta$ ia a primitive root,
  we obtain  $kb_1(n_1 - n_2) \equiv 0 \mod \varphi(q-1)$.  Hence,
\begin{equation}
\chi_1 \leq T \left(1 + \left[\frac{kb_1T}{\varphi(q-1)} \right] \right). \label{lab7}
\end{equation}
\\Now consider $\chi_2$. By (\ref{lab4}), (\ref{aa7}), $\xi$ is equal to the number of solution of the following equation:
\begin{equation}
\sum_{i=0}^{s-1} \sum_{j=1}^{m} h_{ij} \beta^{kb_1i + j -1}  h_{ij} = 0.    \label{l1}
\end{equation}
Bearing in mind that $(h_{11}, \ldots, h_{sm}) \neq \mathbf{0}$, and (\ref{l1}) is a polynomial equation on the field, we get: $\xi \leq k b_1s$.
Thus 
\begin{equation}
\chi_2 \leq \frac{k b_1s T^2}{\varphi(q-1)}.  \label{lab5}
\end{equation}
By (\ref{aa16}), (\ref{lab7}), (\ref{lab5}) and Lemma 6, we have
\[
\chi^2 \leq T \left(1 + \left[\frac{kb_1T}{\varphi(q-1)}\right]\right) + \frac{k b_1s T^2}{\varphi(q-1)} \leq 
     T \left( 1+ (s+1) \frac{k b_1 T}{\varphi(q-1)}\right)
\]
$$
 \leq T(1 + (s+1)b_1 40 \ln3 T  \ln(3 \ln 3T)) \leq  40 (s+2)b_1 T \ln 3T \ln(3 \ln 3T).
$$
We see that
\begin{equation}
\frac{2}{\pi} m \ln p + \frac{7}{5} m - \frac{m-1}{p} \leq 2.5 \ln T, \quad {\rm for} \quad m =[\log_p T] \geq 1.   \label{i3}
\end{equation} 
From (\ref{5a}), (\ref{5c}), (\ref{aa1}), (\ref{aa2}), Lemma 1 and Lemma3, we obtain:
\[
\frac{1}{q \varphi(q-1)} \sum_{\alpha \in F_q} \sum_{\beta \in \hat{F}_q} T \tilde{D}_T(\alpha, \beta) \leq 
\left(2.5 \ln T \right)^s \left(3 + 2 \ln T\right) 40^{1/2}(s+2)^{1/2}  b_1^{1/2}  T^{1/2}
\] 
\begin{equation}
\times   \ln^{1/2} 3T  \ln^{1/2}( 3\ln 3T)  \leq 3 \cdot 40^{1/2} (s+2)^{1/2} 2.5^s   b_1^{1/2}  T^{1/2} \ln^{s+1.5} 3T  \ln^{1/2} (3\ln3 T).  \label{lab8}
\end{equation}
\\Let $T_i =4^i,\;\;\; i \geq [\log_4p] +1$, and let
\begin{equation}
R(\alpha, \beta) = \sum_{b_1 =1}^q\sum_{b_2 =1}^q \sum_{i=[\log_4p] +1}^{\left[\log_4 4q\right]}
 \frac{(81 b_1^{1.5} b_2 \ln^2 3b_1 \ln^2 3b_2)^{-1}  T_i^{1/2} \tilde{D}_{T_i}(\alpha, \beta)}
 { 40^{1/2} (s+2)^{1/2} 2.5^s   \ln^{s+2.5} 3T_i  \ln^{2.5} (3\ln 3T_i)  }.   \label{lab9}
\end{equation}
By (\ref{lab8}), we have:
$$
  \frac{1}{q \cdot \varphi(q-1)} \sum_{\beta \in \hat{F}_q} \sum_{\alpha \in F_q} R(\alpha, \beta)
$$
\begin{equation}
 \leq    \sum_{b_1 =1}^q\sum_{b_2 =1}^q \sum_{i=1}^{\infty}
 \frac{1}
 {9 b_1 b_2 \ln^2 3b_1 \ln^2 3b_2 \; 3i \ln 4  \;\ln^2 (i \ln4)}<  1. \label{lab10}
\end{equation}
Let
\begin{equation}
\Omega_\epsilon = \left\{\alpha \in F_q, \beta \in \hat{F}_q \; \Big| \; R(\mathbf{\alpha}, \mathbf{\beta}) < 1/\epsilon \right\},\;\;\;\;\;
\#\Omega_\epsilon = \gamma q\cdot \varphi(q-1).  \label{aa10}
\end{equation}
Let's prove, that 
\begin{equation}
\gamma \geq 1 - \epsilon.         \label{lab11}   
\end{equation}
We see that $q\cdot \varphi(q-1) (1- \gamma)$  is the number of $\alpha \in F_q, \beta \in \hat{F}_q$, such that $R(\alpha, \beta) \geq \frac{1}{\epsilon}$.
From (\ref{lab10})  we obtain
\[
  1 \geq \frac{1}{q \cdot \varphi(q-1)} \sum_{\alpha \in F_q} \sum_{\beta \in \hat{F}_q} R(\alpha, \beta) \geq \frac{1}{q \cdot \varphi(q-1)} \sum_{(\alpha, \beta) \in \Omega_\epsilon^c} R(\alpha, \beta) 
\] 
\[
\geq \frac{1}{q \cdot \varphi(q-1)} \sum_{(\alpha, \beta) \in \Omega_\epsilon^c} \frac{1}{\epsilon}
= \frac{1}{q \cdot \varphi(q-1)} \frac{1}{\epsilon} \#\Omega_\epsilon^c 
\]
\[
=\frac{1}{q \cdot \varphi(q-1)} \frac{1}{\epsilon} (1 - \gamma) q \cdot \varphi(q-1) = \frac{1 - \gamma}{\epsilon}.
\]
The inequality (\ref{lab11}) is proved.

Now, let $N \in [T_{i_0}, T_{i_0+1})$ for some  $i_0 \in [\log_4p , \log_4 q]$, where $T_i =4^i$, $q>3000$. From (\ref{aa}), (\ref{lab9}) and (\ref{lab10}) we have for all $(\alpha, \beta) \in \Omega_{\epsilon}$,
\[
N D(\mathbf{x}(n, \alpha, \beta))_{0 \leq n \leq {N-1}} \leq s p + T_{i_0+1} \tilde{D}_{T_{{i_0+1}}}(\alpha, \beta) \leq 
s p +\epsilon^{-1}  3^4 40^{1/2} (s+2)^{1/2}  2.5^s 
\]
\begin{equation} \label{300}
\times   T_{i_0+1}^{1/2} \ln^{s+2.5} 3T_{i_0+1}   \ln^{2.5} (3\ln 3T_{i_0+1})  \leq \epsilon^{-1}  c_1 
   N^{1/2} \ln^{s+2.5} 6N  \ln^{2.5}(3 \ln 6N) ,  
\end{equation}
with
$$
    c_1 = s p +   3^6  (s+2)^{1/2} 2.5^s      b_1^{1.5}b_2\ln^2 3b_1 \ln^2 3b_2.
$$
It is easy to see that $c_1 \geq \max(p,3000)$. Hence,  if  $N \leq \max (p,3000)$, then (\ref{300}) is also true. Theorem 1 is proved. \qed
\section{Proof of the Theorem 2}
By (\ref{w1}) and Theorem A, with $m = [\log_p N_1 \ldots N_r]$ and $s_0=s_1 \cdots s_r$, we get:
\begin{equation}
N_1  \cdots  N_r  D ((\mathbf{x}_{\mathbf{n}})_{1 \leq n_w < N_w,\;  1 \leq w \leq r}) \leq s_0 p  
+ \sum_{H \in C^*(p)^{s_0 \times m}} W_p(H) \left|S(H)  \right|,  \nonumber
\end{equation}
where
\[
S(H) = \sum_{ \substack{ n_i \in [0, N_i-1] \\ i=1,...,r}}   e\left(\frac{1}{p}
 \sum_{ \substack{ i_w \in [0, s_w-1] \\ w=1,...,r}} 
  \sum_{j=1}^{k} \sum_{l=1}^{r} h_{i_1, \ldots, i_r, j}  Tr \left(\alpha_l \beta_l^{k (n_l + i_l + s \widetilde{n_l+i_l})+j-1} \right)\right),
\]
\\Let $N_i \in [1, T_i],\;\;\; T_i \in [1, q], \;\;\; 1 \leq i \leq r$.
Using Lemma 4, we get:
\[
S(H)
\leq \sum_{m_{1} = -T_1/2}^{T_1/2} \ldots \sum_{m_{r} = -T_r/2}^{T_r/2} \frac{\left|\dot{S}(H, T, \boldsymbol{\alpha}, \boldsymbol{\beta}) \right|}{\bar{m}_{1} \ldots \bar{m}_{r}},
\]
where
\begin{equation}
\dot{S}(H, T, \boldsymbol{\alpha}, \boldsymbol{\beta}) = \left|\sum_{t_1=0}^{T_1-1} \ldots \sum_{t_r=0}^{T_r-1} e\left(\sigma(\mathbf{t}, \boldsymbol{\alpha}, \boldsymbol{\beta}) \right)  \right|,    \label{aa14}
\end{equation}
with
\[
\sigma(\mathbf{t}, \boldsymbol{\alpha}, \boldsymbol{\beta}) = \sum_{l=1}^{r} \Big( \frac{1}{p} 
 \sum_{i_1 = 0}^{s_1 -1} \ldots \sum_{i_r = 0}^{s_r-1} \sum_{j=1}^{k} h_{i_1, \ldots, i_r, j}  Tr \left(\alpha_l \beta_l^{k (t_l + i_l + s \widetilde{t_l + i_l})+j-1}\right)+ \frac{m_{l} t_l}{T_l}\Big).
\]
Hence, 
\begin{equation}
N_1 \ldots N_r \cdot D((\mathbf{x}_{\mathbf{n}})_{0 \leq n_w < N_w,\;  1 \leq w \leq r}) \leq s_0 p + T_1 \ldots T_r \widetilde{D}_{T_1, \ldots, T_r}(\boldsymbol{\alpha}, \boldsymbol{\beta}),    \label{aa15}
\end{equation}
where,
\begin{equation}
T_1 \ldots T_r \widetilde{D}_{T_1, \ldots, T_r}(\boldsymbol{\alpha}, \boldsymbol{\beta}) = \sum_{H \in C^*(p)^{s_0 \times m}} W_p(H) \sum_{m_{1} = -T_1/2}^{T_1/2} \ldots \sum_{m_{r} = -T_r/2}^{T_r/2} \frac{\left|\dot{S}(H, T, \boldsymbol{\alpha}, \boldsymbol{\beta}) \right|}{\bar{m}_{1} \ldots \bar{m}_{r}},  \nonumber
\end{equation}
with $m = [\log_p (T_1 \ldots T_r)]$. 
Let
\begin{equation} \label{600}
   \chi(\mathbf{T}) = \frac{1}{q^r (\varphi(q-1))^r} \sum_{\boldsymbol{\beta} \in \hat{F}_q^r} 
\sum_{\boldsymbol{\alpha} \in F_q^r}  \left|\dot{S}(H, T, \boldsymbol{\alpha}, \boldsymbol{\beta}) \right|.
\end{equation}
\\Using the Cauchy - Shwartz inequality, we get:
\begin{equation}
\chi^2(\mathbf{T}) \leq \frac{1}{q^r \cdot \left(\varphi(q-1)\right)^r}
 \sum_{\boldsymbol{\beta} \in \hat{F}_q^r} 
\sum_{\boldsymbol{\alpha} \in F_q^r} 
 \left|\sum_{t_1=0}^{T_1-1} \ldots \sum_{t_r=0}^{T_r-1} e\left(\sigma(\mathbf{t}, \boldsymbol{\alpha}, \boldsymbol{\beta}) \right)  \right|^2. \label{aa12}
\end{equation}
By (\ref{aa14}), we have 
\[
\chi^2(\mathbf{T}) = \frac{1}{q^r \cdot \left(\varphi(q-1)\right)^r} \sum_{\boldsymbol{\beta} \in \hat{F}_q^r} \sum_{\boldsymbol{\alpha} \in F_q^r} \sum_{t_1^{(1)},t_1^{(2)} = 0}^{T_1 - 1} \ldots \sum_{t_r^{(1)},t_r^{(2)} = 0}^{T_r - 1} 1
\]
\[
\times e\left(\frac{1}{p} 
\sum_{i_1 = 0}^{s_1-1} \ldots \sum_{i_r = 0}^{s_r-1} \sum_{j=1}^{k} \sum_{l=1}^{r}h_{i_1, \ldots, i_r, j} \left( Tr \left(\alpha_l \beta_l^{k(t_l^{(2)}+i_l + s \widetilde{t_l^{(2)}+i_l})+j-1}\right) \right) \right)
\]
\[
\times e\left(-  Tr \left(\alpha_l \beta_l^{k(t_l^{(1)}+i_l + s \widetilde{t_l^{(1)}+i_l})+j-1}\right) + \frac{m_{1} (t_1^{(2)}-t_1^{(1)})}{T_1} + \ldots + \frac{m_{r} (t_r^{(2)}- t_r^{(1)})}{T_r}\right). 
\]
Using Lemma 2, we get
\[
\chi^2(\mathbf{T})\leq \frac{1}{\left(\varphi(q-1)\right)^r} \sum_{\boldsymbol{\beta} \in \hat{F}_q^r} \sum_{t_1^{(1)},t_1^{(2)} = 0}^{T_1 - 1} \ldots \sum_{t_r^{(1)},t_r^{(2)} = 0}^{T_r - 1} 
\]
\[
\prod_{l=1}^{r} \delta \Big(\sum_{i_1 = 0}^{s_1-1} \ldots \sum_{i_r = 0}^{s_r-1} \sum_{j=1}^{k} h_{i_1, \ldots, i_r, j} \big(\beta_l^{k(t_l^{(2)}+i_l + s \widetilde{t_l^{(2)}+i_l})+j-1} -  \beta_l^{k(t_l^{(1)}+i_l + s \widetilde{t_l^{(1)}+i_l})+j-1} \big) \Big)
\]
\[
\times e\left(\frac{m_{1} (t_1^{(2)}-t_1^{(1)})}{T_1} + \ldots +\frac{m_{r} (t_r^{(2)}- t_r^{(1)})}{T_r}\right)
\]
\[
\leq \frac{1}{\left(\varphi(q-1)\right)^r} \sum_{\boldsymbol{\beta} \in \hat{F}_q^r} \sum_{t_1^{(1)},t_1^{(2)} = 0}^{T_1 - 1} \ldots \sum_{t_r^{(1)},t_r^{(2)} = 0}^{T_r - 1} 
\]
\begin{equation}
\prod_{l=1}^{r} \delta \Big(\sum_{i_1 = 0}^{s_1-1} \ldots \sum_{i_r = 0}^{s_r-1} \sum_{j=1}^{k} h_{i_1, \ldots, i_r, j} \big(\beta_l^{k(t_l^{(2)}+i_l + s \widetilde{t_l^{(2)}+i_l})+j-1} -  \beta_l^{k(t_l^{(1)}+i_l + s \widetilde{t_l^{(1)}+i_l})+j-1} \big) \Big). \nonumber
\end{equation}
It is easy to see  that
\[
\chi^2(\mathbf{T}) \leq \frac{1}{(\varphi(q-1))^r} \sum_{\substack{ \beta_i \in \hat{F}_q \\ i = 1, \ldots, r }} \sum_{\substack{ t_i^{(j)} \in [0, T_i -1]\\ i = 1, \ldots, r, \;\;\; j=1,2}} \prod_{l=1}^{r} \delta (\varsigma_l),
\]
with
\[
\varsigma_l = \sum_{i_1 = 0}^{s_1-1} \ldots \sum_{i_r = 0}^{s_r-1} \sum_{j=1}^{k} h_{i_1, \ldots, i_r, j} \left(\beta_l^{k(t_l^{(2)} + s \widetilde{t_l^{(2)}+i_l})+j-1} -  \beta_l^{k(t_l^{(1)}+i_l + s \widetilde{t_l^{(1)}+i_l})+j-1} \right).
\]
We take a new variable $v_l$ instead of $\widetilde{t_l^{(2)}+i_l},\;\;(l = 1, \ldots, r)$. Enlarging the domain of the summation, we obtain:
\[
\chi^2(\mathbf{T}) \leq \frac{1}{(\varphi(q-1))^r} \sum_{\substack{ \beta_i \in \hat{F}_q \\ i = 1, \ldots, r }} \sum_{\substack{ t_i^{(j)} \in [0, T_i -1]\\ i = 1, \ldots, r, \;\;\; j=1,2}} \prod_{l=1}^{r} \sum_{v_l \in [0, s^{r-1})} \delta (\zeta_l),
\]
where
\begin{equation}
\zeta_l = \sum_{i_1 = 0}^{s_1-1} \ldots \sum_{i_r = 0}^{s_r-1} \sum_{j=1}^{k} h_{i_1, \ldots, i_r, j} \left(\beta_l^{k(t_l^{(2)} + s v_l + i_l)+j-1} -  \beta_l^{k(t_l^{(1)}+i_l + s \widetilde{t_l^{(1)}+i_l})+j-1} \right).  \label{aa161}
\end{equation}
Hence,
\begin{equation}
\chi^2(\mathbf{T}) \leq \prod_{l=1}^{r} \sum_{t_i^{(1)} \in [0, T_i -1]} \xi_l(\mathbf{t}^{(1)}), \label{aa17}
\end{equation}
with
\begin{equation}
\xi_l(\mathbf{t}^{(1)}) = \frac{1}{\varphi(q-1)} \sum_{\beta_l \in \hat{F}_q} \sum_{t_l^{(2)} \in [0, T_l -1]} \sum_{v_l \in [0, s^{r-1})} \delta (\zeta_l).    \label{ll2}
\end{equation}
Consider the equation $\zeta_l = 0$. By (\ref{aa161}), we have:
\begin{equation}
\beta_l^{k((t_l^{(2)}- t_l^{(1)})} \gamma_1 = \gamma_2, \label{l1a}
\end{equation}
where
\[
\gamma_1 = \sum_{i_1 = 0}^{s_1 -1} \ldots \sum_{i_r = 0}^{s_r-1} \sum_{j=1}^{k} h_{i_1, \ldots, i_r, j} \beta_l^{k(s v_l + i_l) + j-1},
\]
\begin{equation}
\gamma_2 = \gamma_2(\beta_l) =  \sum_{i_1 = 0}^{s_1 -1} \ldots \sum_{i_r = 0}^{s_r-1} \sum_{j=1}^{k} h_{i_1, \ldots, i_r, j} \beta_l^{k(s \widetilde{t_l^{(1)}+i_l} + i_l)+j-1}.   \label{aa16a}
\end{equation}
Similarly to (\ref{aa16}), we get:
\begin{equation}
\xi_l(\mathbf{t}^{(1)}) \leq \dot{\xi_l}(\mathbf{t}^{(1)})+ \ddot{\xi_l}(\mathbf{t}^{(1)}),  \label{aa17a}
\end{equation}
where
\[
\dot{\xi_l}(\mathbf{t}^{(1)}) = \frac{1}{\varphi(q-1)} \sum_{\beta_l \in \hat{F}_q} \sum_{t_l^{(2)} \in [0, T_l -1]} \sum_{v_l \in [0, s^{r-1})} 
\delta \left( \beta_l^{k((t_l^{(2)}- t_l^{(1)})} \gamma_1 - \gamma_2\right) 
\]
\begin{equation}
\times \left(1 - \delta(\gamma_2) \right),\;\;\;\;\;\rm{and}\;\;\;\;\;\ddot{\xi_l}(\mathbf{t}^{(1)}) = T_l \frac{1}{\varphi(q-1)} \sum_{\beta_l \in \hat{F}_q} \sum_{v_l \in [0, s^{r-1})} \delta (\gamma_2).
\label{aa19}
\end{equation}
Consider $\dot{\xi_l}(\mathbf{t}^{(1)})$. We see that if $\gamma_2 = 0$, then $\dot{\xi_l}(\mathbf{t}^{(1)}) = 0$. By (\ref{l1a}) and (\ref{aa19}) if $\gamma_2 \neq 0$ and $\gamma_1 = 0$, then also $\dot{\xi_l}(\mathbf{t}^{(1)}) = 0$. Now let $\gamma_1 \neq 0$ and $\gamma_2 \neq 0$. There exists 
 an  integer $a$, such that $\beta_l^a = \gamma_2 / \gamma_1$. By (\ref{l1a}), we get $k (t_l^{(2)} - t_l^{(1)})\equiv$$ a \mod \varphi(q-1)$.
Hence
\[
\# \left\{0 \leq t_l^{(2)} < T_l \;\Big|\; k (t_l^{(2)} - t_l^{(1)}) \equiv a \mod \varphi(q-1) \right\} \leq 1 + \left[\frac{kT_l}{\varphi(q-1)} \right].
\]
By (\ref{aa19}) and Lemma 6, we get
\begin{equation}
\dot{\xi_l}  \leq s^{r-1} (1+    40  \ln3T_l \ln 3(\ln 3T_l)). \label{aa20}
\end{equation}
Similarly to (\ref{lab5}), we have
\[
\# \left\{\beta_l \in \hat{F}_q \Big| \gamma_2(\beta_l)=0 \right\} \leq ks^r, \;\;\;\;\; {\rm and} \;\;\;\;\;\ddot{\xi_l}(\mathbf{t}^{(1)}) \leq 
  k  s^{2r-1}  \frac{T_l}{\varphi(q-1)}.
\]
By (\ref{aa17}), (\ref{aa20}) and Lemma 6 we obtain
\[
\xi_l(\mathbf{t}^{(1)}) \leq  s^{r-1} +    40s^{r-1} \ln3T_l  \ln (3\ln 3T_l)
    +    40s^{2r-1} \ln3T_l  \ln (3 \ln 3T_l)
\]
$$
             \leq    81s^{2r-1} \ln3T_l \ln (3\ln 3T_l).
$$
From (\ref{aa17}), we have
\[
\chi^2(\mathbf{T}) \leq \prod_{l=1}^{r} T_l   81s^{2r-1} \ln3T_l \ln (3\ln 3T_l) .
\]
Using (\ref{i3}), (\ref{600}),  (\ref{aa12}), Lemma1 and Lemma4, we deduce
\[
\frac{1}{q^r (\varphi(q-1))^r} \sum_{\boldsymbol{\alpha} \in F_q^r} \sum_{\boldsymbol{\beta} \in \hat{F}_q^r} T_1 \ldots T_r \tilde{D}_{T_1, \ldots, T_r} (\boldsymbol{\alpha}, \boldsymbol{\beta}) \leq \left(2.5 \ln (T_1 \ldots T_r) \right)^s 
\]
\begin{equation}
\times  \prod_{l=1}^{r} \left(2 + 3 \ln T_l\right) T_l^{1/2} 
    (81s^{2r-1} \ln3T_l \ln (3\ln 3T_l))^{1/2}   \quad {\rm for} \quad T_1\cdots T_r \geq p.  \label{lab12}
\end{equation}
Let $T_{j_i} = 4^{j_i}, \; j_i = 0,1,..., \; i = 1,...,r$, and let
\begin{equation}
R(\boldsymbol{\alpha}, \boldsymbol{\beta})=
     \sum_{\substack{ 1 \leq s_1, \ldots, s_r \leq s \\s_1  \cdots  s_r \leq s }}
 \sum_{\substack{ 1 \leq j_1, \ldots, j_r \leq \log_4 q  \\ \log_4 p \leq j_1 + \ldots + j_r \leq \log_4 4q}}
\frac{ \psi^{-1}(T_\mathbf{j})\sqrt{T_{j_1} \ldots T_{j_r}} \tilde{D}_{T_\mathbf{j}} (\boldsymbol{\alpha}, \boldsymbol{\beta})}{
    3^{5r}s_1 \cdots s_r \ln^2 3s_1 \cdots \ln^2 3s_r}, \label{aaa} 
\end{equation}
where
\[
\psi(\mathbf{T}_\mathbf{j}) =  \left(2.5 \ln (T_1 \ldots T_r) \right)^s \prod_{l=1}^r  s^{r-1/2}
 \ln^{2.5} 3T_{j_l}  \ln^{2.5} (3\ln 3T_{j_l}).
\]
From (\ref{lab12}) we get 
$$
\frac{1}{ q^r \left(\varphi(q-1)\right)^r} \sum_{\mathbf{\alpha} \in F_q^r} \sum_{\mathbf{\beta} \in \hat{F}_q^r} R(\mathbf{\alpha}, \mathbf{\beta}) 
$$
\begin{equation}
 \leq   \prod_{i=1}^r  \sum_{s_i =1}^{\infty}  \sum_{j_i=1}^{\infty}
 \frac{1}
 {3s_i \ln^2 3s_i \cdot  3j_i \ln 4  \;\ln^2 (j_i \ln4)}<  1. \nonumber
\end{equation}
Let
\begin{equation}
\Omega_\epsilon = \left\{\boldsymbol{\alpha} \in F_q^r, \boldsymbol{\beta} \in \hat{F}_q^r \; \Big| \; R(\boldsymbol{\alpha}, \boldsymbol{\beta}) < \frac{1}{\epsilon} \right\}, \;\;\;\;\; \#\Omega_\epsilon = \gamma q^r \left(\varphi(q-1)\right)^r.   \label{lab15}
\end{equation}
Similarly to (\ref{lab10})-(\ref{lab11}), we get $\gamma \geq 1 - \epsilon$.

Now, let $N_i \in [T_{j_i}, T_{j_i+1})$ for some $j_i \in [0, \log_4 q], \; i=1,...,r$, with
 $T_{j_i} =4^{j_i}$ and  $T_{j_1} \ldots T_{j_r} \geq p$. From (\ref{aa15}), (\ref{aaa}) and (\ref{lab15}), we have for all  $(\boldsymbol{\alpha}, \boldsymbol{\beta}) \in \Omega_{\epsilon}$:
$$
N_1 \ldots N_r D \left((\mathbf{x}_{\mathbf{n}})_{1 \leq n_w < N_w,\;  1 \leq w \leq r} \right) 
  \leq s_0 p  + T_{j_1 +1} \ldots T_{j_r+1} \tilde{D}_{T_\mathbf{j}} (\boldsymbol{\alpha}, \boldsymbol{\beta})
$$
$$
\leq s_0 p + \epsilon^{-1}  3^{5r}  s^{r^2 - r/2} \left(T_{j_1+1} \ldots T_{j_r+1} \right)^{1/2} 
$$
\[
\times \left(2.5 \ln \left(T_{j_1+1} \ldots T_{j_r+1} \right) \right)^s
\prod_{l=1}^{r}   \ln^{2.5} (3T_{j_l+1})  \ln^{2.5} (3\ln 3T_{j_l+1} )  s_l \ln^2 3s_l 
\]
\begin{equation} \label{302}
\leq  \epsilon^{-1}c_2 \left(N_1 \ldots N_r \right)^{1/2} 
 \ln^{s+2.5r} (2^{r+1} N_1 \ldots N_r)     \ln^{2.5r} 3\ln (6N_1 \cdots N_r)
\end{equation}
with
$$
  c_2 = s_0 p +  3^{5r}  2.5^s  s^{r^2 - r/2}   s_1 \cdots s_r \ln^2 3s_1 \cdots \ln^2 3s_r 
$$

It is easy to see that $c_2 \geq \max(p,3000) $. Hence, if $N_1 \cdots N_r \leq \max (p,3000)$, then (\ref{302}) is also true. Theorem 2 is proved. \qed

Address: Department of Mathematics,
Bar-Ilan University,
Ramat-Gan, 52900 Israel \\
e-mails: $ mlevin@math.biu.ac.il,$\;\;\;$ volinskaya\_i@yahoo.com$.

\end{document}